\documentclass[a4paper]{article}
\usepackage{color}
\usepackage{amsmath}
\usepackage{amssymb}
\newtheorem{theorem}{Theorem}
\newtheorem{prop}{Proposition}
\newtheorem{defi}{Definition}

\newcommand{\ba}{\begin{eqnarray}}
\newcommand{\ea}{\end{eqnarray}}
\newcommand{\ban}{\begin{eqnarray*}}
\newcommand{\ean}{\end{eqnarray*}}
\newcommand{\no}{\nonumber}

\newcommand{\sg}{\sigma}

\newcommand{\mapright}[1]{%
\smash{\mathop{%
\hbox to 1.0cm{\rightarrowfill}}\limits^{#1}}}
\newcommand{\mapleft}[1]{%
\smash{\mathop{%
\hbox to 1.3cm{\leftarrowfill}}\limits^{#1}}}

\begin{document}
\title{
\begin{flushright}
  \begin{minipage}[b]{5em}
    \normalsize
    ${}$      \\
  \end{minipage}
\end{flushright}
{\bf Harmonic Partitions of Positive Integers and Bosonic Extension of Euler's Pentagonal Number Theorem }}
\author{Masao Jinzenji${}^{(1)}$, Yu Tajima${}^{(2)}$\\
\\${}^{(1)}$ \it Department of Mathematics,  \\
\it Okayama University \\
\it  Okayama, 700-8530, Japan\\
\\${}^{(2)}$ \it Division of Mathematics, Graduate School of Science \\
\it Hokkaido University \\
\it  Kita-ku, Sapporo, 060-0810, Japan\\
\\
\it e-mail address: (1) pcj70e4e@okayama-u.ac.jp \\
\it\hspace{2.0cm}(2) yu149573@icloud.com}
\maketitle

\begin{abstract}
In this paper, we first propose a cohomological derivation of the celebrated Euler's Pentagonal Number Theorem.
Then we prove an identity that corresponds to a bosonic extension of the theorem. 
The proof corresponds to a cohomological re-derivation of Euler's another celebrated identity. 

\end{abstract}

\section{Introduction}
First part of this paper is a presentation of a naive cohomological approach that one of the authors (M. J.) obtained 
in trying to solve a problem in the book ``Combinatorial Problems and Exercises'' by L. Lov\.{a}sz \cite{lovasz}, that requests readers to 
prove {\bf  Euler's Pentagonal Number Theorem \cite{euler}}:
\ba
\prod_{m=1}^{\infty}(1-q^m)=1+\sum_{l=1}^{\infty}(-1)^{l}\biggl(q^{\frac{l(3l-1)}{2}}+q^{\frac{l(3l+1)}{2}}\biggr).
\label{pent}
\ea
Let $S_{l}(n)$ be set of partitions of positive integer $n$ with $l$ distinct parts.
\ban
S_{l}(n)&:=&\{(n_1,n_2,\cdots ,n_l)\; |\; n_i\in\mathbb{N}\;\;(i=1,2,\cdots ,l), \;\;n_1>n_2>\cdots >n_l>0, \\
&&n_1+n_2+\cdots +n_l =n \}.
\ean
Then the l.h.s. of (\ref{pent}) is represented as follows:
\ba
\prod_{m=1}^{\infty}(1-q^m)=1+\sum_{n=1}^{\infty}\biggl(\sum_{l=0}^{\infty}(-1)^{l}(S_{l}(n))^{\sharp}\biggr)q^{n}.
\ea
M.J. was very interested in elementary cohomology theory in his graduate student days, and he was tempted to interpret
$\sum_{l=0}^{\infty}(-1)^{l}(S_{l}(n))^{\sharp}$ as some Euler number of a complex. Then he defined a real vector 
space $CS_{l}(n)$ spanned by elements of $S_{l}(n)$ and tried to define coboundary operator $\delta: CS_{l}(n)
\rightarrow CS_{l+1}(n)$. This trial turned out to be successful. Moreover, adjoint $\delta^{\ast}: CS_{l}(n)
\rightarrow CS_{l-1}(n)$ of the coboundary operator $\delta$ can be defined, and Laplacian $\delta\delta^{\ast}+
\delta^{\ast}\delta$ on $CS_{l}(n)$ can also be defined. Then we can define harmonic partition of $n$ with distinct parts.
These harmonic partitions exist only if $n=\frac{l(3l-1)}{2}$ or $n=\frac{l(3l+1)}{2}$ and explicitly written as follows:
\ba
n=\frac{l(3l-1)}{2}:&& (2l-1,2l-2,\cdots,l+1,l),\no\\
n=\frac{l(3l+1)}{2}:&& (2l,2l-1,\cdots,l+2,l+1).
\ea  
(\ref{pent}) immediately follows from these constructions. In Section 2, we give detailed explanation of  them. 
For combinatorial proof of Euler's pentagonal number theorem related to our approach, we have to mention  Franklin's classical work \cite{frank}.

With these results, we tried to extend the above naive constructions to the case of ordinary partitions. Let $P_{l}(n)$ be set of ordinary partitions of positive integer $n$ with $l$ parts.
\ba
P_{l}(n):=\{(N_1,N_2,\cdots ,N_{l}) \;|\; N_1\geq N_2\geq\cdots \geq N_{l}>0, \;\sum_{j=1}^{l}N_{j}=n\;\}.
\ea
We immediately have,
\ba
\prod_{m=1}^{\infty}\frac{1}{(1+q^m)}=1+\sum_{n=1}^{\infty}\biggl(\sum_{l=0}^{\infty}(-1)^{l}(P_{l}(n))^{\sharp}\biggr)q^{n}.
\ea
Then we defined a real vector space $CP_{l}(n)$ spanned by elements of $P_{l}(n)$ and tried to extend the above
constructions to $CP_{l}(n)$. This trial turned out to be successful.  We found that harmonic ordinary partitions are explicitly given by, 
\ba
&&((n_1)^{n_1+2l},(n_2)^{2t_2},\cdots ,(n_k)^{2t_k})\no\\
&&(n_1>n_2>\cdots >n_k>0, \;\;l\geq 0,\;\; t_2,\cdots ,t_k\geq1).
\ea
Therefore, we obtain the following identity:
\ba
\prod_{m=1}^{\infty}\frac{1}{(1+q^m)}=1+\sum_{l=1}^{\infty}(-1)^{l}\frac{q^{l^2}}{\displaystyle{\prod_{j=1}^{l}(1-q^{2j})}}.
\label{main}
\ea
Section 3 is devoted to cohomological proof of this identity. 

We have to point out that it is equivalent to the well-known result by Euler \cite{euler}:
\ba
\prod_{m=1}^{\infty}(1+q^{2m-1})=1+\sum_{l=1}^{\infty}\frac{q^{l^2}}{\displaystyle{\prod_{j=1}^{l}(1-q^{2j})}}.
\label{euler2}
\ea
We can derive (\ref{main}) from (\ref{euler2}) by changing $q$ into $-q$ and using the well-known identity by Euler \cite{euler}:
\ba
\prod_{m=1}^{\infty}\frac{1}{(1-q^{2m-1})}=\prod_{m=1}^{\infty}(1+q^m).
\ea  
Therefore, the discussion in this note corresponds to cohomological re-derivation of (\ref{euler2}).

\section{The Case of Partitions with Distinct Parts and Euler's Pentagonal Number Theorem}

\begin{defi}
Let us define the partition of  positive integer $n$ with distinct parts:\\
\ba
&&\sg:=(n_1,n_2,\cdots ,n_l),\no\\
&&(n_i\in\mathbb{N}\;\;(i=1,2,\cdots ,l),\;\; n_1>n_2>\cdots >n_l>0, n_1+n_2+\cdots +n_l =n).\no
\ea
We call $l$ length of the partition $\sg$. Then, we denote by $S_{l}(n)$ set of partitions of $n$ with $l$ distinct parts:
\ban
S_{l}(n)&:=&\{(n_1,n_2,\cdots ,n_l)\; |\; n_i\in\mathbb{N}\;\;(i=1,2,\cdots ,l), \;\;n_1>n_2>\cdots >n_l>0, \\
&&n_1+n_2+\cdots +n_l =n \}
\ean
Then we define as $CS_l(n)$ a real vector space generated by elements of $S_l(n)$: 
\ban
CS_l(n):=\bigoplus_{\sg\in S_l(n)}\mathbb{R}\sg.
\ean
\end{defi}

\begin{defi}
Let $\sg=(n_1,\cdots ,n_{l})\in CS_{l}(n)$ be a partition of $n$ with distinct parts. We define as $m(\sg)$ the minimum positive integer that satisfies $n_{m+1}<n_m-1$. In other words, $m(\sg)$ is the maximum positive integer that satisfies,
\ban
n_{j}=n_{j-1}-1,\;\;(j=2,3,\cdots,m(\sg)).
\ean 
If $l=1$, we define $m(\sg)=1$. If $n_{j}=n_{j-1}-1$ for $j=2,3,\cdots,l$, we define $m(\sg)=l$.
\end{defi}

\begin{defi}
We define a linear map $\delta : CS_{l}(n)\rightarrow CS_{l+1}(n)$ that 
is determined by action on each generator $\sg=(n_1,\cdots ,n_{l})\in CS_{l}(n)$:
\begin{itemize}
\item[{(i)}]The case of $m(\sg)<{l}$\\
$m(\sg)\geq n_l\Rightarrow \delta (\sg)=0.$\\
$m(\sg)<n_{l}\Rightarrow \delta (\sg)=(n_1-1,n_2-1,\cdots ,n_{m(\sg)}-1,n_{m(\sg)+1},\cdots ,n_{l},m(\sg))$.
\item[{(ii)}]The case of $m(\sg)=l$\\
$m(\sg)\geq n_{l}-1\Rightarrow \delta(\sg)=0.$\\
$m(\sg)<n_{l}-1\Rightarrow \delta(\sg)=(n_1-1,n_2-1,\cdots ,n_{m(\sg)}-1,m(\sg))$.
\end{itemize}
\end{defi}

\begin{prop}
\label{coho1}
$\delta\circ\delta=0$.
\end{prop}
{\bf Proof)} Let $\sg\in CS_l(n)$ be a partition of $n$ with distinct parts. If $\delta(\sg)=0$, then $\delta(\delta(\sg))=0$ obviously holds. Let us assume that $\delta(\sg)\neq 0$ and that $\delta(\sg)$ is written as $(N_1,\cdots ,N_{l+1})$. Since $m(\delta(\sg))\geq m(\sg)=N_{l+1}$, $\delta (\delta (\sg))=0$ holds by the definition of $\delta$. $\Box$

\begin{defi}
 We define a linear map $\delta^{\ast} : CS_{l}(n)\rightarrow CS_{l-1}(n)$ 
that 
is determined by action on each generator $\sg=(n_1,\cdots ,n_{l})\in CS_{l}(n)$:
\begin{itemize}
\item[{(i)}]The case of $m(\sg)<{l}$\\
$m(\sg)< n_l\Rightarrow \delta^{\ast} (\sg)=0$.\\
$m(\sg) \geq n_{l}\Rightarrow \delta^{\ast} (\sg)=(n_1+1,n_2+1,\cdots ,n_{n_{l}}+1,n_{n_{l}+1},\cdots ,n_{{l}-1})$.
\item[{(ii)}]The case of $m(\sg)=l$\\
$m(\sg)\leq n_{l}\Rightarrow \delta^{\ast} (\sg)=0$.\\
$m(\sg)>n_{l}\Rightarrow \delta^{\ast}(\sg)=(n_1+1,n_2+1,\cdots ,n_{n_{l}}+1,n_{n_{l}+1}\cdots ,n_{l-1})$.
\end{itemize}
\end{defi}

\begin{prop}
\label{adj1-1}
$\delta^{\ast}\circ\delta^{\ast}=0$.
\end{prop}
{\bf Proof)} Let $\sg\in CS_l(n)$ be a partition of $n$. If $\delta^{\ast}(\sg)=0$, then $\delta^{\ast}(\delta^{\ast}(\sg))=0$ obviously holds. Let us assume that $\delta^{\ast}(\sg)\neq 0$ and that $\delta^{\ast}(\sg)$ is written as $(N_1,\cdots ,N_{l-1})$. Since $m(\delta^{\ast}(\sg))=n_{l}< n_{l-1}=N_{l-1}$, $\delta^{\ast} (\delta^{\ast} (\sg))=0$ holds by the definition of $\delta^{\ast}$. $\Box$

\begin{prop}
\label{adj1-2}
\quad
\begin{itemize}
\item[{(i)}]{$\delta(\sg)=\tau\neq 0 \Rightarrow \delta^{\ast}(\tau)=\sg$}
\item[{(ii)}]{$\delta^{\ast}(\tau)=\sg\neq 0 \Rightarrow \delta(\sg)=\tau$}
\end{itemize}
\end{prop}
{\bf Proof)} 
\begin{itemize}
\item[{(i)}]Let $\sg=(n_1,n_2,\cdots ,n_l)\in CS_{l}(n)$ be a partition of $n$ with distinct parts.  We assume that $\delta(\sg)\neq 0$, and $\delta (\sg)=(N_1,N_2,\cdots ,N_{l+1})=\tau$. Then, $m(\tau)\geq m(\sg)=N_{l+1}$.\\
(a) The case of $m(\tau)<l+1$\\
$\delta^{\ast}(\tau)=(N_1+1,N_2+1,\cdots ,N_{N_{l+1}}+1,N_{N_{l+1}+1},\cdots ,N_l)\\
=(N_1+1,N_2+1,\cdots ,N_{m(\sg)}+1,N_{m(\sg)+1},\cdots ,N_l)\\
=(n_1,n_2,\cdots ,n_{m(\sg)},n_{m(\sg)+1},\cdots ,n_l)=\sg$.\\
(b) The case of $m(\tau)=l+1$\\
Since $l+1>m(\sg)$, $m(\tau)>m(\sg)=N_{l+1}$. Therefore, $\delta^{\ast}(\tau)$ equals $\sg$ in the same way as tha case (a).

\item[{(ii)}] Let $\tau=(n_1,n_2,\cdots ,n_l)\in CS_{l}(n)$ be a partition of $n$ with distinct parts.  We assume that $\delta^{\ast}(\tau)\neq 0$, and $\delta^{\ast}(\tau)=(N_1,N_2,\cdots ,N_{l-1})=\sg$. Then, $m(\sg)= n_{l}<n_{l-1}$.\\
(a) The case of $m(\sg)<l-1$\\
Since $m(\sg)<n_{l-1}=N_{l-1}$,\\
$\delta(\sg)=(N_1-1,N_2-1,\cdots ,N_{m(\sg)}-1,N_{m(\sg)+1},\cdots ,N_{l-1},m(\sg))\\
=(N_1-1,N_2-1,\cdots ,N_{n_l}-1,N_{n_{l}+1},\cdots ,N_{l-1},n_l)\\=(n_1,n_2,\cdots ,n_{n_l},n_{n_{l}+1},\cdots ,n_{l-1},n_l)=\tau$.\\
(b) The case of $m(\sg)=l-1$\\
Since $m(\sg)<n_{l-1}=N_{l-1}-1$, $\delta(\sg)$ equals $\tau$ in the same way as the case (a). $\Box$
\end{itemize}

\begin{prop}
\label{hm1}
Let $\sg=(n_1,n_2,\cdots,n_l)\in CS_l(n)$ be a partition of $n$ with distinct parts.
Then we have,
\ban
&&\delta(\sg)=\delta^{\ast}(\sg)=0\\
&\Longleftrightarrow& m(\sg)=l, n_l=m(\sg),\;\; \mbox{or}\;\; m(\sg)=l, n_l=m(\sg)+1.
\label{harmd}
\ean
\end{prop}
{\bf Proof)}\\
$(\Longleftarrow)$ We assume that $m(\sg)=l$.\\
In the case of $n_l=m(\sg)$, $\delta(\sg)=\delta^{\ast}(\sg)=0$ holds because $m(\sg)=n_l>n_l-1$.\\
In the case of $n_l=m(\sg)+1$, $\delta(\sg)=\delta^{\ast}(\sg)=0$ holds because $m(\sg)=n_l-1<n_l$.\\
$(\Longrightarrow)$ Let $\sg=(n_1,n_2,\cdots,n_l)\in CS_l(n)$ be a partition of $n$ with distinct parts. We assume that $\delta(\sg)=\delta^{\ast}(\sg)=0$. From the condition $\delta(\sg)=0$, we are naturally led to consider the following two cases.
\begin{itemize}
\item[{(i)}]{In the case of $n_l \leq m(\sg)<l$, $\delta(\sg)$ vanishes.  However, we conclude that $\delta^{\ast}(\sg)\neq 0$ by the definition of $\delta^{\ast}$.}
\item[{(ii)}]{In the case of $n_l-1\leq m(\sg)=l$, $\delta(\sg)$ vanishes. The condition $\delta^{\ast}(\sg)=0$ tells us that $m(\sg)\leq n_l$. Since $n_l-1\leq m(\sg)\leq n_l$, $m(\sg)=n_l-1$ or $n_l$.}
\end{itemize}
Therefore, $m(\sg)=l, n_l=m(\sg),\;\; \mbox{or}\;\; m(\sg)=l, n_l=m(\sg)+1$.  $\Box$\\
\begin{defi}
From Proposition \ref{coho1}, we can define cohomology $H^{l}(CS(n))$ as 
$\mbox{Ker}(\delta|_{CS_{l}(n)})/\mbox{Im}(\delta|_{CS_{l-1}(n)})$. Then,
\ban
\sum_{l:even}(S_{l}(n))^{\sharp}-\sum_{l:odd}(S_{l}(n))^{\sharp},
\ean
can be interpreted as Euler number:
\ba
\chi(H^{*}(CS(n))):=\sum_{l:even}\dim_{\bf R}(H^{l}(CS(n)))-\sum_{l:odd}\dim_{\bf R}(H^{l}(CS(n))). 
\ea
We define symmetric positive definite inner product of $CS_{l}(n)$ by,
\ba
(\sigma,\tau)=\delta_{\sigma\tau},\;\;(\sigma,\tau\in S_{l}(n)).
\label{inner1}
\ea
Then Propositions \ref{adj1-1} and \ref{adj1-2} tell us,
\ba
(\delta\sigma,\tau)=(\sigma,\delta^{\ast}\tau),
\label{adjoint1}
\ea
i.e., $\delta$ and $\delta^{*}$ are adjoint to each other.   
Then we define a Laplacian $\delta\delta^{\ast}+\delta^{\ast}\delta$ on $CS_{l}(n)$.
Since we have (\ref{inner1}) and (\ref{adjoint1}), we can immediately conclude,
\ba
(\delta\delta^{\ast}+\delta^{\ast}\delta)\sigma=0\Longleftrightarrow \delta(\sigma)=\delta^{\ast}(\sigma)=0.
\ea 
We define a linear space of harmonic partitions ${\bf H}^{l}(CS(n))$ by $\mbox{Ker}((\delta\delta^{\ast}+\delta^{\ast}\delta)|_{CS_{l}(n)})$. Since we have standard decomposition,
\ba
CS_{l}(n)={\bf H}^{l}(CS(n))\oplus \mbox{Im}(\delta|_{CS_{l-1}(n)})\oplus \mbox{Im}(\delta^{\ast}|_{CS_{l+1}(n)}),
\ea
we obtain,
\ba
\dim_{\bf R}({\bf H}^{l}(CS(n)))=\dim_{\bf R}(H^{l}(CS(n))).
\ea
We call $\sg\in CS_l(n)$ that satisfies $\delta(\sg)=\delta^{\ast}(\sg)=0$,
which is nothing but a linear base of ${\bf H}^{l}(CS(n))$, a harmonic partition of $n$ with
distinct parts. Proposition \ref{hm1} tells us that harmonic partitions with distinct parts are explicitly given as follows.  
\ba
n=\frac{l(3l-1)}{2}:&& (2l-1,2l-2,\cdots,l+1,l),\no\\
n=\frac{l(3l+1)}{2}:&& (2l,2l-1,\cdots,l+2,l+1).
\ea 
\end{defi}

\begin{theorem}({\bf  Euler's Pentagonal Number Theorem})
\ba
\prod_{m=1}^{\infty}(1-q^m)=1+\sum_{l=1}^{\infty}(-1)^{l}\biggl(q^{\frac{l(3l-1)}{2}}+q^{\frac{l(3l+1)}{2}}\biggr).
\ea
\end{theorem}
{\bf Proof)} 
\ban
\prod_{m=1}^{\infty}(1-q^m)&=&1+\sum_{n=1}^{\infty}\biggl(\sum_{l:even}(S_{l}(n))^{\sharp}-\sum_{l:odd}(S_{l}(n))^{\sharp}\biggr)q^{n}\\
&=&1+\sum_{n=1}^{\infty}\biggl(\sum_{l:even}\dim_{\bf R}(H^{l}(CS(n)))-\sum_{l:odd}\dim_{\bf R}(H^{l}(CS(n)))\biggr)q^{n}\\
&=&1+\sum_{n=1}^{\infty}\biggl(\sum_{l:even}\dim_{\bf R}({\bf H}^{l}(CS(n)))-\sum_{l:odd}\dim_{\bf R}({\bf H}^{l}(CS(n)))\biggr)q^{n}\\
&=&1+\sum_{l=1}^{\infty}(-1)^{l}\biggl(q^{\frac{l(3l-1)}{2}}+q^{\frac{l(3l+1)}{2}}\biggr). \;\;\Box
\ean

\section{The Case of Ordinary Partitions and Bosonic Extension of Euler's Pentagonal Number Theorem } 

\begin{defi}
Let us define the ordinary partition of  positive integer $n$ as follows :
\ban
&&\sg=(N_1,N_2,\cdots ,N_{l})\\
&&(N_i\in\mathbb{N}\;\;(i=1,2,\cdots ,l),\; N_1\geq N_2\geq\cdots \geq N_{l}>0,\; \sum_{j=1}^{l}N_{j}=n\;).
\ean
We call $l$ length of the partition $\sg$. Then, we denote by $P_{l}(n)$ set of the ordinary partitions of $n$ with length $l$ :
\ban
P_{l}(n):=\{(N_1,N_2,\cdots ,N_{l}) \;|\; N_1\geq N_2\geq\cdots \geq N_{l}>0, \;\sum_{j=1}^{l}N_{j}=n\;\}.
\ean
We define as $CP_l(n)$ a real vector space generated by elements of $P_l(n)$ : \\
$CP_l(n):=\bigoplus_{\sg\in P_l(n)}\mathbb{R}\sg$.
\end{defi}

\begin{defi}
Let us represent a partition $\sg$ of the positive integer n as follows:
\ban
&&\sg=(n_1^{m_1},n_2^{m_2},\cdots,n_k^{m_k})\no\\ 
&& (n_1>n_2>\cdots>n_k>0,\;\;n=\sum_{j=1}^{k}n_j m_j, \;l=\sum_{j=1}^{k}m_j\;).
\ean
We call each $(n_t)^{m_t}(t=1,2,\cdots,k)\ $a block.
 Then, we define a map $\delta:CP_l(n)\rightarrow CP_{l+1}(n)$ as follows: \\
First, we devide the block $(n_1)^{m_1}$ into $(n_1-1)^{m_1}$ and $(1)^{m_1}$ and transpose $(1)^{m_1}$ into $(m_1)^1$
(from now on, we simply call this operation dividing  $(n_1)^{m_1}$ into  $(n_1-1)^{m_1}$ and $(m_1)^1$).
\begin{itemize}
\item[{(i)}] { If $m_1>n_1-1$, then  $\delta (\sg)=0$. } 
\item[{(ii)}] { If there exists $j\;(2\leq j\leq k)$ such that $m_1=n_j$ and if $m_j,m_{j+1},\cdots,m_k$ are all even, then we put the block  $(m_1)^1$ on the right side of $((n_j)^{m_j}):$
\ban
\delta(\sg)=((n_1-1)^{m_1},(n_2)^{m_2},\cdots,(n_j)^{m_j+1},(n_{j+1})^{m_{j+1}},\cdots,(n_k)^{m_k}).
\ean
 If there exists some  odd number in $m_j, m_{j+1},\cdots,m_k$, then $\delta(\sg)=0$.}
\item[\mbox{(iii)}] { In the case except for the cases of  (i), (ii),\\
 there exists unique integer $i\;\;(1\leq i\leq k)$ that satisfies $n_i>m_1>n_{i+1}$.\\
 If $m_{i+1},m_{i+2},\cdots,m_k$ are all even, or $i=k\;\;(n_k>m_1)$, then we put the block $(m_1)^1$ on the right side of $(n_i)^{m_i}:$
\ban
\delta(\sg)=((n_1-1)^{m_1},(n_2)^{m_2},\cdots,(n_i)^{m_i},(m_1)^1,(n_{i+1})^{m_{i+1}},\cdots,(n_k)^{m_k}).
\ean
 If there exists some odd number in $m_{i+1}, m_{i+2},\cdots,m_k$, then $\delta(\sg)=0$.} 
\end{itemize}
Let us remark a subtle point.\\
If there appear in $\delta(\sg)=(N_1^{M_1},N_2^{M_2},\cdots,N_K^{M_K})$ blocks $(N_T)^{M_T},(N_{T+1})^{M_{T+1}}\\
(T=1,2,\cdots,K-1)$ that satisfy $N_T=N_{T+1}$, we automatically rewrite\\
$\delta(\sg)= (N_1^{M_1},\cdots,N_T^{M_T},N_{T+1}^{M_{T+1}},\cdots,N_K^{M_K})$ into $(N_1^{M_1},\cdots,N_T^{M_T+M_{T+1}},\cdots,N_K^{M_K})$.
\end{defi}

\begin{prop}
\label{coho2}
$\delta\circ\delta=0$.
\end{prop}
{\bf Proof)} Let $\sg\in CP_l(n)$ be a partition of $n$. If $\delta(\sg)=0$, then $\delta(\delta(\sg))=0$ obviously holds. Let us assume that $\delta(\sg)\neq 0$ and that $\delta(\sg)$ is written as $(N_1^{M_1},\cdots,N_J^{M_J},N_{J+1}^{M_{J+1}},\cdots,N_K^{M_K})$. 
Then, we have the following two cases.
\begin{itemize}
\item[(i)] {There exists odd number $M_J$ with $M_{J+1},\cdots,M_K$ all even, and $M_1\geq N_J$. }
\item[(ii)] { $M_1,\cdots, M_{K}$ are all even and $M_{1}>N_{1}-1$.}
\end{itemize}
In both cases, $\delta(\delta(\sg))=0$ holds by definition of $\delta$. $\Box$

\begin{defi}
Let $\sg=(n_1^{m_1},n_2^{m_2},\cdots,n_k^{m_k})\;\;\; (n_1>n_2>\cdots>n_k>0,\,\;n=n_1 m_1+n_2 m_2+\cdots+n_k m_k)$
be a partition in $CP_l(n)$.  Then we define a map $\delta^{\ast} :CP_l(n)\rightarrow CP_{l-1}(n)$ 
in the following way.\\
First, if there exist some odd numbers among $m_j$'s $(1\leq j\leq k)$, we denote by $t$ the maximum index $j$ of odd $m_{j}$'s. \\
(i) If $t=1$ or $m_1,m_2,\cdots,m_k$ are all even, then we divide the block $(n_1)^{m_1}$ into $(n_1)^{m_1-1}$ and $(n_{1})^{1}$ and transpose $(n_{1})^{1}$ into $(1)^{n_1}$ (from now on, we simply call this operation dividing  $(n_1)^{m_1}$ into $(n_1)^{m_1-1}$ and 
$(1)^{n_1}$). 
\begin{itemize}
\item[{(a)}]{ In the case of $n_1>m_1$ or $m_1=n_1+2i\;\;(i\geq 0)$, we define $\delta^{\ast}(\sg)=0$.}
\item[{(b)}]{In the case of $m_1=n_1+2i+1\;\;(i\geq 0)$, we put the block of $(1)^{n_1}$ left-aligned on $(n_1)^{m_1-1}$,
\ban 
\delta^{\ast}(\sg)=(({n_1}+1)^{n_1},(n_1)^{2i},(n_2)^{m_2},\cdots ,(n_k)^{m_k}).
\ean
}
\end{itemize}
(ii) If $t\geq 2$, then we divide the block $(n_t)^{m_t}$ into $(n_t)^{m_t-1}$ and $(1)^{n_t}$.
\begin{itemize}
\item[{(a)}]{In the case of $n_t>m_1$, we define $\delta^{\ast}(\sg)=0$.}
\item[{(b)}]{In the case of $n_t\leq {m_1}$, we put the block $(1)^{n_t}$ left-aligned on $(n_1)^{m_1}$, 
\ban
\delta^{\ast}(\sg)&=&(({n_1}+1)^{n_t},(n_1)^{m_1-n_t},(n_2)^{m_2},\cdots\\
&&\cdots,(n_{t-1})^{m_{t-1}},(n_t)^{m_t-1},(n_{t+1})^{m_{t+1}},\cdots ,(n_k)^{m_k}).
\ean 
}
\end{itemize}
Let us remark that if there appear a block $(n_t)^0$, then we automatically omit it.
\end{defi}

\begin{prop}
\label{adj2-1}
$\delta^{\ast}\circ\delta^{\ast}=0$.
\end{prop}
{\bf Proof)} Let $\sg\in CP_l(n)$ be a partition of $n$. If $\delta^{\ast}(\sg)=0$, then $\delta^{\ast}(\delta^{\ast}(\sg))=0$ obviously holds. Let us assume that $\delta^{\ast}(\sg)\neq 0$.
Then, if $\delta^{\ast}(\sg)$ is written as $(N_1^{M_1},\cdots,N_T^{M_T},N_{T+1}^{M_{T+1}},\cdots,N_K^{M_K})$, 
we can easily see that there exists odd number $M_T$ with $M_{T+1},\cdots,M_K$ all even number and $N_T> M_1$, or that $M_1, \cdots ,M_{K}$ are all even and $N_1>M_1$. Therefore, by definition of $\delta^{\ast}$,  $\delta^{\ast}(\delta^{\ast}(\sg))=0$ also holds. $\Box$

\begin{prop}
\label{adj2-2}
\quad
\begin{itemize}
\item[{(i)}]{$\delta(\sg)=\tau\neq 0 \Rightarrow \delta^{\ast}(\tau)=\sg$}
\item[{(ii)}]{$\delta^{\ast}(\tau)=\sg\neq 0 \Rightarrow \delta(\sg)=\tau$}
\end{itemize}

\end{prop}
{\bf Proof)} 
\begin{itemize}
\item[{(i)}]{Let $\sg=(n_1^{m_1},n_2^{m_2},\cdots,n_k^{m_k})\in CP_l(n)$ be a partition of $n$. We assume that $\delta(\sg)=\tau \neq 0$ and $\delta(\sg)=(N_1^{M_1},N_2^{M_2},\cdots ,N_K^{M_K})=\tau$.
We denote by $(N_{J+1})^{M_{J+1}}$ the block on the right side of the block where $(m_1)^1$ is added by $\delta$. Then $N_J=m_1$, $M_J$ is an odd number and $M_{J+1},\cdots ,M_{K}$ are all even. 
Let us  apply $\delta^{\ast}$ to $\tau$. Then we divide the block of $(N_J)^{M_J}$ into $(N_J)^{M_J-1}$ and $(1)^{N_J}=(1)^{m_1}$, and put the block $(1)^{N_J}$ left-aligned on $(N_1)^{M_1}$. This operation corresponds to putting the block $(m_1)^1$, that was moved by $\delta$, back to its original position. Therefore, $\delta^{\ast}(\tau)=\sg$.}

\item[{(ii)}]{Let $\tau=(n_1^{m_1},n_2^{m_2},\cdots,n_k^{m_k})\in CP_l(n)$ be a partition of $n$. We assume that $\delta^{\ast}(\tau)=\sg \neq 0$ and $\delta^{\ast}(\tau)=(N_1^{M_1},N_2^{M_2},\cdots ,N_K^{M_K})=\sg$.
By the map $\delta^{\ast}$, we put the block of $(1)^{n_t}$ left-aligned on $(n_1)^{\tilde{m_1}}$. (Note that ${\tilde{m_1}}=m_1-1$ if $t=1$, and ${\tilde{m_1}}=m_1$ if $t\neq 1$.) Then, $M_1=n_t$. The blocks $(n_{t+1})^{m_{t+1}},\cdots ,(n_k)^{m_k}$ remain the same, and we represent them as $(N_{J+1})^{M_{J+1}},\cdots ,(N_K)^{M_K}$.  Then $M_{J+1},\cdots ,M_K$ are all even. Next, by the map $\delta$, we divide the block $(N_1)^{M_1}$ into $(N_1-1)^{M_1}$ and $(M_1)^1=(n_t)^1$ and put the block of $(n_t)^1$ on the left side of $(N_{J+1})^{M_{J+1}}$. This operation corresponds to putting the block of $(1)^{n_t}$, which was moved by $\delta^{\ast}$,  back to its original position. Therefore, $\delta(\sg)=\tau$.}

\end{itemize}
$\Box$

\begin{prop}
\label{hm2}
Let $\sg=(n_1^{m_1},n_2^{m_2},\cdots,n_k^{m_k})\in CP_l(n)$ be a partition of $n$.
Then we have,
\ba
&&\delta(\sg)=\delta^{\ast}(\sg)=0\no\\
&\Longleftrightarrow& \sg=((n_1)^{n_1+2l},(n_2)^{2t_2},\cdots ,(n_k)^{2t_k})\no\\
&&(n_1>n_2>\cdots >n_k>0, \;\;l\geq 0,\;\; t_2,\cdots ,t_k\geq 1).
\label{harm}
\ea
\end{prop}
{\bf Proof)}\\
$(\Longleftarrow)$ Let $\sg=((n_1)^{n_1+2l},(n_2)^{2t_2},\cdots ,(n_k)^{2t_k})\;\;(n_1>n_2>\cdots >n_k, \;\;l\geq 0, \;\;t_2,\cdots ,t_k\geq 1)=:(N_1^{M_1},\cdots,N_K^{M_K})$ be a partition of $n$. Since $n_1+2l=M_1>N_1-1=n_1-1$, $\delta(\sg)$ vanishes. 
On the other hand, since $M_1=n_1+2l=N_1+2l\;\;(l\geq0)$ and $M_{2},\cdots, M_{K}$ are all even, $\delta^{\ast}(\sg)$ also 
vanishes.\\
$(\Longrightarrow)$ Let $\sg=(n_1^{m_1},n_2^{m_2},\cdots,n_k^{m_k})\in CP_l(n)$ be a partition of $n$. We assume that $\delta(\sg)=\delta^{\ast}(\sg)=0$. From the condition $\delta(\sg)=0$, we are naturally led to consider the following three
cases.
\begin{itemize}
\item[{(i)}]{In the case of $m_1>n_1-1$, $\delta(\sg)$ vanishes.  Since $n_1\leq m_1$, the condition $\delta^{\ast}(\sg)=0$
tells us that $m_2,\cdots ,m_k$ are all even and that $m_1=n_1+2l\;\;(l\geq 0)$. Therefore, $\sg$
is represented as $((n_1)^{n_1+2l},(n_2)^{2t_2},\cdots ,(n_k)^{2t_k})$.}

\item[{(ii)}]{If there exist $j \;\;(2\leq j\leq k)$ that satisfy $m_1=n_j$ and some odd numbers among $m_j,\cdots ,m_k$, then $\delta(\sg)$ vanishes. Let $m_t$ be the odd number that have maximum index $t$. Since $2\leq j\leq t\leq k$ and $n_t\leq m_1=n_j$, 
we conclude that $\delta^{\ast}(\sigma)\neq 0$.}

\item[{(iii)}]{If there exist unique integer $i\;\;(1\leq i < k)$ that satisfies $n_i>m_1>n_{i+1}$ and some odd numbers among $m_{i+1},\cdots ,m_k$, then $\delta(\sg)$ vanishes. Let $m_t$ be the odd number that have maximum index $t$. Since $2\leq i+1\leq t\leq k$ and $n_i>m_1>n_{i+1}\geq n_t$, we conclude that $\delta^{\ast}(\sigma)\neq 0$.}
\end{itemize}
Therefore, $\sigma$ takes the form presented in the last line of (\ref{harm}).  $\Box$

\begin{defi}
From Proposition \ref{coho2}, we can define cohomology $H^{l}(CP(n))$ as 
$\mbox{Ker}(\delta|_{CP_{l}(n)})/\mbox{Im}(\delta|_{CP_{l-1}(n)})$. Then,
\ban
\sum_{l:even}(P_{l}(n))^{\sharp}-\sum_{l:odd}(P_{l}(n))^{\sharp},
\ean
can be interpreted as Euler number:
\ba
\chi(H^{*}(CP(n))):=\sum_{l:even}\dim_{\bf R}(H^{l}(CP(n)))-\sum_{l:odd}\dim_{\bf R}(H^{l}(CP(n))). 
\ea
We define symmetric positive definite inner product of $CP_{l}(n)$ by,
\ba
(\sigma,\tau)=\delta_{\sigma\tau},\;\;(\sigma,\tau\in P_{l}(n)).
\label{inner2}
\ea
Then Propositions \ref{adj2-1} and \ref{adj2-2} tell us,
\ba
(\delta\sigma,\tau)=(\sigma,\delta^{\ast}\tau),
\label{adjoint2}
\ea
i.e., $\delta$ and $\delta^{*}$ are adjoint to each other.   
Then we define a Laplacian $\delta\delta^{\ast}+\delta^{\ast}\delta$ on $CP_{l}(n)$.
Since we have (\ref{inner2}) and (\ref{adjoint2}), we can immediately conclude,
\ba
(\delta\delta^{\ast}+\delta^{\ast}\delta)\sigma=0\Longleftrightarrow \delta(\sigma)=\delta^{\ast}(\sigma)=0.
\ea 
We define a linear space of harmonic partitions ${\bf H}^{l}(CP(n))$ by $\mbox{Ker}((\delta\delta^{\ast}+\delta^{\ast}\delta)|_{CP_{l}(n)})$. Since we have standard decomposition,
\ba
CP_{l}(n)={\bf H}^{l}(CP(n))\oplus \mbox{Im}(\delta|_{CP_{l-1}(n)})\oplus \mbox{Im}(\delta^{\ast}|_{CP_{l+1}(n)}),
\ea
we obtain,
\ba
\dim_{\bf R}({\bf H}^{l}(CP(n)))=\dim_{\bf R}(H^{l}(CP(n))).
\ea
We call $\sg\in CP_l(n)$ that satisfies $\delta(\sg)=\delta^{\ast}(\sg)=0$,
which is nothing but a linear base of ${\bf H}^{l}(CP(n))$, a harmonic ordinary partition of $n$. 
 According to Proposition \ref{hm2}, harmonic ordinary partitions up to $n=26$ are given as follows.
\ban
n=1:&& (1),\no\\
n=2:&& \mbox{none},\no\\
n=3:&& (1^3),\no\\
n=4:&& (2^2),\no\\
n=5:&& (1^5),\no\\
n=6:&& (2^2,1^4),\no\\
n=7:&& (1^7),\no\\
n=8:&& (2^4),(2^2,1^4),\no\\
n=9:&& (3^3),(1^9)\no\\
n=10:&& (2^4,1^2),(2^2,1^6)\no\\
n=11:&& (3^3,1^2),(1^{11}),\no\\
n=12:&& (2^6),(2^4,1^4),(2^2,1^8),\no\\
n=13:&& (3^3,2^2),(3^3,1^4),(1^{13}),\no\\
n=14:&& (2^6,1^2),(2^4,1^6),(2^2,1^{10}),\no\\
n=15:&& (3^5),(3^3,2^2,1^2),(3^3,1^6),(1^{15}),\no\\
n=16:&& (4^4),(2^8),(2^6,1^4),(2^4,1^8),(2^2,1^{12}),\no\\
n=17:&& (3^5,1^2),(3^3,2^4),(3^3,2^2,1^4),(3^3,1^8),(1^{17}),\no\\
n=18:&& (4^4,1^2),(2^8,1^2),(2^6,1^6),(2^4,1^{10}),(2^2,1^{14}),\no\\
n=19:&& (3^5,2^2),(3^5,1^4),(3^3,2^4,1^2),(3^3,2^2,1^6),(3^3,1^{10}),(1^{19}),\no\\
n=20:&& (4^4,2^2),(4^4,1^4),(2^{10}),(2^8,1^4),(2^6,1^8),(2^4,1^{12}),(2^2,1^{16}),\no\\
n=21:&&(3^7),(3^5,2^2,1^2),(3^5,1^6),(3^3,2^6),(3^3,2^4,1^4),(3^3,2^2,1^8),(3^3,1^{12}),(1^{21}),\no\\
\ean
\ban
n=22:&&(4^4,3^2),(4^4,2^2,1^2),(4^4,1^6),(2^{10},1^2),(2^8,1^6),(2^6,1^10),(2^4,1^{14}),(2^2,1^{18}),\no\\
n=23:&&(3^7,1^2),(3^5,2^4),(3^5,2^2,1^4),(3^5,1^8),(3^3,2^6,1^2),(3^3,2^4,1^6),(3^3,2^2,1^{10}),(3^3,1^{14}),\no\\
&&(1^{23}),\no\\
n=24:&&(4^6),(4^4,3^2,1^2),(4^4,2^4),(4^4,2^2,1^4),(4^4,1^8),(2^{12}),(2^{10},1^4),(2^8,1^8),(2^6,1^{12}),\no\\
&&(2^4,1^{16}),(2^2,1^{20}),\no\\
n=25:&&(5^5), (3^7,2^2), (3^7,1^4),(3^5,2^4,1^2),(3^5,2^2,1^6),(3^5,1^{10}),(3^3.2^8),(3^3,2^6,1^4),\no\\
&&(3^3,2^4,1^8),(3^3,2^2,1^{12}),(3^3,1^{16}),(1^{25}),\no\\
n=26:&&(4^6,1^2),(4^4,3^2,2^2),(4^4,3^2,1^4),(4^4,2^4,1^2),(4^4,2^2,1^6),(4^4,1^{10}),(2^{12},1^2),\no\\
&&(2^{10},1^6),(2^8,1^{10}),(2^6,1^{14}),(2^4,1^{18}),(2^2,1^{22}).
\ean 
\end{defi}

\begin{theorem}
\ba
\prod_{m=1}^{\infty}\frac{1}{(1+q^m)}=1+\sum_{l=1}^{\infty}(-1)^{l}\frac{q^{l^2}}{\displaystyle{\prod_{j=1}^{l}(1-q^{2j})}}.
\ea
\end{theorem}
{\bf Proof)}
\ban
\prod_{m=1}^{\infty}\frac{1}{(1+q^m)}&=&1+\sum_{n=1}^{\infty}\biggl(\sum_{l:even}(P_{l}(n))^{\sharp}-\sum_{l:odd}(P_{l}(n))^{\sharp}\biggr)q^{n}\\
&=&1+\sum_{n=1}^{\infty}\biggl(\sum_{l:even}\dim_{\bf R}(H^{l}(CP(n)))-\sum_{l:odd}\dim_{\bf R}(H^{l}(CP(n)))\biggr)q^{n}\\
&=&1+\sum_{n=1}^{\infty}\biggl(\sum_{l:even}\dim_{\bf R}({\bf H}^{l}(CP(n)))-\sum_{l:odd}\dim_{\bf R}({\bf H}^{l}(CP(n)))\biggr)q^{n}\\
&=&1+\sum_{l=1}^{\infty}(-1)^{l}\frac{q^{l^2}}{\displaystyle{\prod_{j=1}^{l}(1-q^{2j})}}.\;\;\Box
\ean

\end{document}